\documentclass{birkjour}
\usepackage{amssymb,enumerate,array,multirow}
 \newtheorem{thm}{Theorem}[section]
 
 \newtheorem{lem}[thm]{Lemma}
 \newtheorem{prop}[thm]{Proposition}
 \theoremstyle{definition}
 \newtheorem{defn}[thm]{Definition}
 \theoremstyle{remark}
 \newtheorem{rem}[thm]{Remark}
 \newtheorem*{ex}{Example}
 \numberwithin{equation}{section}
 
 \newcommand{\kk}{{\Bbbk}}

\begin{document}

%
%
%
%
%
%
%
%
%

\title[Graded Clifford Algebras and Graded Skew Clifford Algebras]{Graded Clifford Algebras and Graded Skew Clifford Algebras and their Role in the Classification of Artin-Schelter Regular Algebras}

\author[P Veerapen]{Padmini Veerapen}

\address{
Department of Mathematics, P.O.~Box 19408\\
Tennessee Technological University,
Cookeville, 
TN 76019-0408\\}

\email{pveerapen@tntech.edu}



\date{January 1, 2004}

\begin{abstract}
This paper is a survey of work done on $\mathbb{Z}$-graded Clifford algebras (GCAs) and $\mathbb{Z}$-graded \textit{skew} Clifford algebras (GSCAs) \cite{VVW, SV, CaV, NVZ, VVe1, VVe2}. In particular, we discuss the hypotheses necessary for these algebras to be Artin Schelter-regular \cite{AS, ATV1} and show how certain `points' called, point modules, can be associated to them. We may view an AS-regular algebra as a noncommutative analog of the polynomial ring. We begin our survey with a fundamental result in \cite{VVW} that is essential to subsequent results discussed here: the connection between point modules and rank-two quadrics. Using, in part, this connection the authors in \cite{SV} provide a method to construct GCAs with finitely many distinct isomorphism classes of point modules. In \cite{CaV}, Cassidy and Vancliff introduce a quantized analog of a GCA, called a graded \textit{skew} Clifford algebra and Nafari et al. \cite{NVZ} show that most Artin Schelter-regular algebras of global dimension three are either twists of graded skew Clifford algebras of global dimension three or Ore extensions of graded Clifford algebras of global dimension two. Vancliff et al. \cite{VVe1, VVe2} go a step further and generalize the result of \cite{VVW}, between point modules and rank-two quadrics, by showing that point modules over GSCAs are determined by (noncommutative) quadrics of $\mu$-rank at most two. 
\end{abstract}

\maketitle
\section{Introduction} 
In the 1980's, with the appearance of many noncommutative algebras from quantum physics  and a desire to find a noncommutative algebraic geometry that would be as successful as commutative algebraic geometry had been for commutative algebra, M. Artin and W. Schelter \cite{AS} introduced the notion of a regular algebra. A few years later Artin, Tate and Van den Bergh \cite{ATV1, ATV2} classified the generic classes of regular algebras of global dimension three. The main idea behind this classification, introduced by Artin, Tate and Van den Bergh in \cite{ATV1, ATV2}, involved using certain graded modules in place of geometric data, for example, point modules in place of certain points and line modules in place of certain lines. In particular, Artin, Tate, and Van den Bergh showed that such algebras could be associated to certain subschemes $E$ (typically of dimension one) of $\mathbb{P}^2$ where points in the scheme $E$ parametrize certain modules over these algebras called point modules. The technique involved the definition of a quantum analog of the projective plane $\mathbb{P}^2$. Since then progress has been made on the classification of regular algebras of global dimension four. However, this is still an open problem and as it stands now, \textit{quadratic} regular algebras of global dimension four are still unclassified. This is where graded Clifford algebras (GCAs) and graded skew Clifford algebras (GSCAs) come into play.

In the papers we survey, we highlight how GCAs \cite{Lb} and GSCAs \cite{CaV} provide a trove of examples of regular algebras of global dimension three and four with the goal of informing the classification of quadratic regular algebras of global dimension four. We begin with a result where Vancliff et al. \cite{VVW} prove a fundamental connection between GCAs and quadrics of rank at most two, stated precisely in Theorem \ref{vvw}. Using this result, we show how Stephenson et al. \cite{SV} construct examples of regular GCAs with a finite number of point modules. In order to further inform the classification of quadratic Artin-Schelter regular (AS-regular) algebras of global dimension $n$, Cassidy and Vancliff in \cite{CaV} introduced a quantized analog of a graded Clifford algebra called a graded \textit{skew} Clifford algebra (GSCA). 

To determine how useful GSCAs are, Nafari et al. \cite{NVZ} prove that almost all quadratic regular algebras of global dimension three can be classified using GSCAs. Indeed, based on \cite{AS} classification of regular algebras of global dimension three as being of types $A, B, E, H, S_1, S_1'$, and $S_2$, Nafari et al. showed that most quadratic AS-regular algebras of global dimension three are either twists of GSCAs of global dimension three or Ore extensions of graded Clifford algebras of global dimension two. In fact, only AS-regular algebras of type $E$ and an open subset of those of type $A$ are possibly not associated to GSCAs in that way. 

In \cite{VVe2}, Vancliff and the author generalize the connection between GCAs and quadrics of rank at most two from \cite{VVW} to GSCAs. This is achieved by, first, defining a notion of rank, called $\mu$-rank, on noncommutative quadratic forms and second, by showing that point modules over GSCAs can be determined by noncommutative quadratic forms of $\mu$-rank at most two \cite{VVe1}. 

We conclude this survey paper with work being done in \cite{ChV, TV} to construct two families of potentially generic regular algebras of global dimension four where the main objective is to determine the line scheme of these families with the goal of ultimately classifying regular quadratic algebras of global dimension four. 

We note that `classical' Clifford algebras as defined by Lounesto \cite{L} differ from the graded Clifford and graded skew Clifford algebras discussed here. One of the main differences is that the former are $\mathbb{Z}_2$-graded while the latter are $\mathbb{Z}$-graded. Chen and Kang \cite{CK} seek to address the challenges posed by the differences in grading by associating Clifford theory to certain graded spaces.

In Section~\ref{sec2} of this paper, we present definitions of notions critical to the other sections and survey results about graded Clifford algebras (GCAs) from \cite{VVW} and \cite{SV}. In Section~\ref{sec3}, we provide a review of results of \cite{NVZ} and discuss how GSCAs of global dimension three are related to the types of AS-regular algebras classified by \cite{ATV1, ATV2}. We also survey existing results from \cite{VVe1} regarding the $\mu$-rank of a noncommutative quadratic form as a generalization of the notion of rank of quadratic forms in the commutative setting. We discuss results from \cite{VVe2} relating the $\mu$-rank of noncommutative quadratic forms to point modules over graded skew Clifford algebras. 

\bigskip
\bigskip



\section{Graded Clifford Algebras}\label{sec2}

Throughout the article, $\Bbbk$~denotes an algebraically closed field such that char$(\Bbbk)\neq~2$, and $M(n,\ \Bbbk)$ denotes the vector space of $n \times n$ matrices with entries in $\kk$. For a graded $\kk$-algebra~$B$, the span of the homogeneous elements in $B$ of degree $i$ will be denoted $B_i$, and the notation $T(V)$ will denote the tensor algebra on the vector space~$V$.  If $C$ is any ring or vector space, then $C^\times$ will denote the nonzero elements in $C$. We use $R$ to denote the polynomial ring on degree-one generators $x_1, \ldots , x_n$.

\medskip

In this section, we survey results on graded Clifford algebras (GCAs). In particular, we discuss methods used to determine point modules over GCAs. We begin with the definition of a GCA and provide other useful definitions.

\medskip

\begin{defn}\label{gca}\cite[\S 4]{Lb}
Let $M_1, \ldots, M_n \in M_n(\Bbbk)$ be symmetric matrices. A graded Clifford algebra is the $\Bbbk$-algebra $C$ on degree-one generators $x_1, \ldots, x_n$ and on degree-two generators $y_1, \ldots, y_n$ with the following defining relations:
\begin{enumerate}[(a)]
\item $x_ix_j + x_jx_i = \displaystyle \sum_{k=1}^n (M_k)_{ij}y_k$ for all $i, j = 1, \ldots, n$;
\item $y_k$ central for all $k = 1, \ldots, n$.
\end{enumerate}
\end{defn}

\medskip

\begin{ex}
Suppose 
$M_1 = \left[\begin{smallmatrix}
2 & \lambda\\
\lambda & 0
\end{smallmatrix} \right]$ and
$M_2 = \left[\begin{smallmatrix}
0 & 0\\
0 & 1
\end{smallmatrix} \right]$, where $\lambda \in \Bbbk$. Let $C$ be the graded $\Bbbk$-algebra on degree-one generators $x_1, x_2$ and on degree-two generators $y_1, y_2$ with defining relations given by Definition \ref{gca}(i). That is, for $i, j = 1, 2$, the defining relations of $C$ are: $2x_1^2 = 2y_1$, $x_2^2 = y_2$, and $x_1x_2 + x_2x_1 = \lambda y_1 = \lambda x_1^2$. Thus, $\frac{\Bbbk\langle x_1, x_2 \rangle}{\langle x_1x_2 + x_2x_1 - \lambda x_1^2 \rangle}$ maps onto $C$. 
\end{ex}

\medskip

\begin{defn}\cite{AS}\label{AS}
Let $A = \oplus_{i=0}^{\infty} A_i$ be a finitely generated, $\mathbb{N}$-graded, connected $\Bbbk$-algebra. The algebra $A$ will be called regular (or AS-regular) if it satisfies the following properties:
\begin{enumerate}
\item[(a)] A has finite global dimension $d$,
\item[(b)] A has polynomial growth, and
\item[(c)] A is Gorenstein.
\end{enumerate}  
\end{defn}

\medskip

\begin{rem}

An equivalent condition to condition (b) is that $A$ has finite Gelfand-Kirillov dimension (see \cite{KL}). Condition (c) is the main reason why an AS-regular algebra is viewed as a noncommutative analogue to the polynomial ring inasmuch as it replaces the symmetry condition of commutativity, as shown in the example below. Henceforth, an AS-regular algebra will be called a regular algebra.
\end{rem}

\medskip

\begin{ex}
Consider the $\Bbbk$-algebra $A = \Bbbk[x, y]$ where $xy = qyx$ and $q \in \Bbbk^{\times}.$ Note that $A$ satisfies the Gorenstein condition. Let $(\ast)$ be the minimal projective resolution of the trivial left $A$-module $_{A}\Bbbk$
\[0 \hspace{1mm} {\xrightarrow{\hspace*{1cm}}}  \hspace{1mm} _{A}A \overset{h}{\xrightarrow{\hspace*{1cm}}}  \hspace{1mm} _{A}A^2 \overset{f}{\xrightarrow{\hspace*{1cm}}}  \hspace{1mm} _{A}A {\xrightarrow{\hspace*{1cm}}}  \hspace{1mm} _{A}\Bbbk {\xrightarrow{\hspace*{1cm}}}  \hspace{1mm} 0 \quad \quad \quad \quad (\ast),\]
where $h(a) = a\left[\begin{matrix}-qy & x\end{matrix} \right]$ and $f(b, c) = \left[\begin{matrix}b & c\end{matrix} \right]\left[\begin{matrix}x\\y\end{matrix} \right]$ for $a, b, c \in A$. This minimal projective resolution is the shortest possible of this form for this $A$. If we now dualize $(\ast)$ by applying the functor $GHom(., _{A}A)$, we obtain the following minimal resolution of the trivial right $A$-module $\Bbbk_A$
\[0 \hspace{1mm} {\xrightarrow{\hspace*{1cm}}} \hspace{1mm} A_A \overset{f^{\ast}}{\xrightarrow{\hspace*{1cm}}}  \hspace{1mm} A_A^2 \overset{h^{\ast}}{\xrightarrow{\hspace*{1cm}}}  \hspace{1mm} A_A {\xrightarrow{\hspace*{1cm}}} \hspace{1mm} \Bbbk_A {\xrightarrow{\hspace*{1cm}}} \hspace{1mm} 0,\]
where $f^{\ast}(c) = \left[\begin{matrix}x\\y\end{matrix} \right]c$ and $h^{\ast} = \left[\begin{matrix}-qy & x\end{matrix} \right]\left[\begin{matrix}a\\b\end{matrix} \right]$ for $a, b, c \in A$. Note that $f^{\ast}$ and $h^{\ast}$ are now left multiplication by the same matrices as in $(\ast)$.
\end{ex}

\medskip

\begin{rem}\label{KDGCA}

\begin{enumerate}[(i)] \item Suppose $B$ is a quadratic regular algebra where $V = B_1$ and $W \subset V \otimes_{\Bbbk} V$ be the span of the defining relations of $B$ and let $\mathcal{V}(W)$ be the zero locus in $\mathbb{P}(V^{\ast}) \times \mathbb{P}(V^{\ast})$ of elements of $W$. The Koszul dual \cite{ST}, $B^{\ast}$, of $B$ is the quadratic algebra given by the quotient of the tensor algebra on $V^{\ast}$ by the ideal generated by $W^{\perp}$ where $W^{\perp} = \{f \in V^{\ast} \otimes_{\Bbbk} V^{\ast}: f(w) = 0 \mbox{ for } w \in W\}$.
\item With $\mathcal{V}(W)$ as in part (i), by \cite{ATV1} $\mathcal{V}(W)$ is the graph of an automorphism $\tau$ of a scheme $E \subset \mathbb{P}^2$.

\item We may connect geometry to a GCA by considering the quadratic form associated to the symmetric matrix $M_k$. This is discussed in \cite{VVW, CaV} where a quadratic form $q_k \in \Bbbk[z_1, \ldots, z_n]$ is associated to each $M_k$ via $q_k = z^TM_kz$ with $z^T = [z_1 \ldots z_n]$. In particular, viewing the basis $\{z_1, \ldots, z_n\}$ as the dual basis to $\{x_1, \ldots, x_n\}$, the (commutative) polynomial algebra $\Bbbk[z_1, \ldots, z_n]$, $R$, is the Koszul dual of the quadratic algebra $B/\langle y_1, \ldots, y_n \rangle$. Moreover, the zero locus $\mathcal{V}(q_k)$ of each $q_k$ is a \textit{quadric} in $\mathbb{P}(C_1^{\ast})$. A \textit{base point} of the quadric system $\mathcal{V}(q_1), \ldots, \mathcal{V}(q_n)$ is a point in $\bigcap_{k=1}^n \mathcal{V}(q_k)$. If no such point exists, then the quadric system is said to be base point free.
\end{enumerate}
\end{rem}

Results in \cite{AL, Lb} relate GCAs to regular algebras using the geometry we discuss above. 

\medskip

\begin{prop}\label{AL&Lb}\cite{AL, Lb}
The graded Clifford algebra $C$ is quadratic, Auslander-regular of global dimension $n$ and satisfies the Cohen-Macaulay property with Hilbert series $1/(1 - t)^n$ if and only if the associated quadric system is base point free; in this case, $C$ is Artin-Schelter regular and is a noetherian domain.
\end{prop}

\medskip

\begin{ex} Let
$M_1 = \left[\begin{smallmatrix}
2 & \lambda\\
\lambda & 0
\end{smallmatrix} \right]$ and
$M_2 = \left[\begin{smallmatrix}
0 & 0\\
0 & 1
\end{smallmatrix} \right]$, where $\lambda \in \Bbbk$. Since the quadratic forms $q_1$ and $q_2$ associated to $M_1$ and $M_2$, respectively, are base point free, then the $\Bbbk$-algebra $\frac{\Bbbk\langle x_1, x_2 \rangle}{\langle x_1x_2 + x_2x_1 - \lambda x_1^2 \rangle}$ is isomorphic to the Clifford algebra, $C$, by Proposition \ref{AL&Lb}.
\end{ex}

\subsection{Graded Clifford Algebras and Quadrics of Rank Two}\label{sec3}

\medskip

We now want to determine certain modules called, point modules, over regular graded Clifford algebras. To that end, we examine Theorem \ref{vvw} in \cite{VVW}. It is helpful to note that in one of their seminal papers \cite{ATV1} prove that, under certain conditions, point modules are determined by a scheme. This scheme represents the functor of point modules and was called in \cite{VVr}, the \textit{point scheme}. The point scheme of regular algebras will be essential to the classification of such algebras in this article.


\begin{defn}\cite{ATV1}
A right (respectively, left) point module over a graded algebra, $B$, is a cyclic graded right (respectively, left) $B$-module $M = \bigoplus_{i \ge 0} M_i$ such that $M = M_0B$ or (respectively, $BM_0$) and $\dim_{\Bbbk}(M_i) = 1$ for all $i$.
\end{defn}

\begin{thm}\label{vvw}\cite[Theorem 1.7]{VVW}
Let $C$ denote a GCA determined by symmetric matrices, $N_1, \ldots, N_n \in M(n, \Bbbk)$ and let $\mathcal{Q}$ be the corresponding (commutative) quadric system in $\mathbb{P}^{n-1}$. If $\mathcal{Q}$ is base point free, then the number of isomorphism classes of left (respectively, right) point modules over $C$ is equal to $2r_2 + r_1 \in \mathbb{N} \cup \{0, \infty\}$ where $r_j$ denotes the number of matrices in $\mathbb{P}\left(\sum_{k=1}^n \Bbbk N_k\right)$ that have rank $j$. If the number of left (respectively, right) point modules is finite, then $r_1 \in \{0, 1\}$. 
\end{thm}

\noindent Theorem \ref{vvw} is critical in that it allows us to control the number of point modules over regular GCAs using the the quadrics of rank at most two in the Koszul dual of the GCA.


\subsection{Graded Clifford Algebras with Finitely Many Points}\label{sec3}

\medskip

We will now discuss how to construct certain regular graded Clifford algebras of global dimension four with a finite number of point modules by reviewing results in \cite{SV}. The authors in \cite{SV} use Theorem \ref{vvw} to construct regular GCAs of global dimension four with a finite number of points in their scheme. Indeed, if a GCA has at least two and at most finitely many distinct isomorphism classes of point modules, then \cite{SV} show that the number of intersection points of two planar cubic divisors determines the number of isomorphism classes of point modules over such GCAs. The two planar cubic divisors, in turn, parametrize quadrics of rank at most two. This is where Theorem \ref{vvw} plays a critical role. For this subsection, we will call a linear system of quadrics $\mathfrak{q}$, such that $\mathbb{P}(\mathfrak{q}) \cong \mathbb{P}^2$, a net of quadrics and a linear system of quadrics $\mathfrak{Q}$, such that $\mathbb{P}(\mathfrak{Q}) \cong \mathbb{P}^3$, a web of quadrics. 

Using Proposition \ref{AL&Lb}, the authors determine whether a GCA of global dimension four is regular. They do so by ensuring that the quadric system associated to the GCA is base point free. The method used to construct regular GCAs with a finite number of points in their point scheme, involves the use of two planar cubic divisors. Suppose $\mathfrak{Q} = \{q_1, q_2, q_3, q_4\}$ where $q_1 = \mathcal{V}(a_1a_2)$ is a web of quadrics in $\mathbb{P}^3$ corresponding to a regular GCA of global dimension four. Recall that from the discussion in Remark \ref{KDGCA}, we can associate a linear system of quadrics to a GCA via its Koszul dual. We consider the quadratic form $Q = a_1a_2 \in R_2$ and the net $\mathfrak{q} = \Bbbk q_2 \oplus \Bbbk q_3 \oplus \Bbbk q_4$ such that $Q \notin \mathfrak{q}$. We view $\mathfrak{q}_i$ as the isomorphic image of the net $\mathfrak{q}$ determined by the map $R_2 \to R_2/R_1a_i$ and let the schemes $F_i$ be the family of quadrics in $\mathbb{P}(\mathfrak{q})$ that meet $\mathcal{V}(a_i)$ in a degenerate conic. In fact, $F_i$ parametrizes the quadrics in $\mathbb{P}(\mathfrak{q}_i)$ of rank two and $\dim(F_i) = 1$ for $i = 1, 2$. The following result states precisely how the intersection points of the cubic divisors determine the number of points in the point scheme of the GCA:

\begin{thm}\cite[Theorem 2.6]{SV}
Let A be a graded Clifford algebra of global dimension four and let $\mathfrak{Q}$ denote the corresponding quadric system that is base-point free. Suppose $\mathfrak{Q} = \Bbbk Q \oplus \mathfrak{q}$, where $\mathfrak{q}$ is a net of quadrics in $\mathfrak{Q}$ and $Q = \mathcal{V}(ab) \notin \mathfrak{q}$ such that $a, b \in R_1$ are linearly independent. Write $F_1$ and $F_2$ for the two planar cubic divisors in $\mathbb{P}(\mathfrak{q})$ that parametrize the quadrics in $\mathbb{P}(\mathfrak{q})$ that meet the two distinct planes of $Q$ in a degenerate conic such that $|F_1 \cap F_2| < \infty$ (that is, $|F_1 \cap F_2| \le 9$). Let $m$ be the number of distinct points in $F_1 \cap F_2$, and let $r_i$ denote the number of distinct quadrics in $\mathbb{P}(\mathfrak{Q})$ of rank $i$, and suppose $r_1 \le 1$. We have
\begin{enumerate}[(a)]
\item $r_2 < \infty$;
\item $A$ has at most finitely many isomorphism classes of point modules;
\item $m - r_1 \le r_2 \le 2(m - r_1) + 1$.
\end{enumerate}
\end{thm}  

\noindent The following result will be helpful in the example below:

\begin{lem}\cite[Lemma 2.3]{SV}
If $Q$ and $Q'$ are distinct quadrics in $\mathbb{P}^3$, then either $\mathbb{P}(\Bbbk Q \oplus \Bbbk Q')$ consists of quadrics of rank at most two or $\mathbb{P}(\Bbbk Q \oplus \Bbbk Q')$ contains at most three quadrics of rank at most two. Moreover, if $\mathbb{P}(\Bbbk Q \oplus \Bbbk Q')$ consists of quadrics of rank at most two, then if $Q = \mathcal{V}(ab)$, where $a, b \in R_1$ are linearly independent, then $Q' = \mathcal{V}(f)$ where $f \in R_1a$ or $f \in R_1b$ or $f \in \Bbbk^{\times} a^2 + \Bbbk ab + \Bbbk^{\times}b^2$.

\end{lem}

\begin{ex}\cite[Example 5.1]{SV}
The goal in this example is to construct a regular GCA with exactly twenty distinct points in its point scheme. First, for the GCA to be regular, one needs a base-point free web $\mathfrak{Q} \in \mathbb{P}^3$ of quadrics such that $\mathbb{P}(\mathfrak{Q})$ contains exactly ten rank-two quadrics. Recall that the ongoing assumption is that the GCA has more than one point in its point scheme so, using the fact that any regular GCA of global dimension four has $r_1 + 2r_2$ point modules, this means that $\mathfrak{Q}$ will contain at least one quadric of rank two.

Suppose $a, b, c, d \in R_1$ are linearly independent and suppose $\mathfrak{q} = \Bbbk(a^2 - b^2) \oplus \Bbbk(a^2 - c^2) \oplus \Bbbk(a^2 - d^2)$. We see that $\mathbb{P}(\mathfrak{q})$ contains six rank-two quadrics and none of rank one. We set $Q = a^2 - e^2$, for some $e \in S_1$, and define $\mathfrak{Q} = \Bbbk Q \oplus \mathfrak{q}$. If $e$ is taken to equal $a + b + c + d$, then $\mathfrak{Q}$ is a base-point free web of quadrics (if char$(\Bbbk) \ne 3$ or 5) such that $\mathbb{P}(\mathfrak{Q})$ contains exactly ten rank-two quadrics. One can then set $a = \pm b$ and obtain two cubic divisors with their six intersection points. The number of rank-two quadrics are obtained by linearly combining the quadrics corresponding to the intersection points with $a^2 - b^2$.
\end{ex}


\section{Graded \textit{Skew} Clifford Algebras}\label{sec3}


In \cite{CaV}, Cassidy and Vancliff introduced a quantized analog of a graded Clifford algebra, known as a graded \textit{skew} Clifford algebra (GSCAs). The introduction of this quantized analog required several notions to be generalized: $\mu$-symmetric matrices, noncommutative quadratic forms, noncommutative quadric system, and base points. We discuss these generalizations and survey papers that generalize results from the previous section, regarding graded Clifford algebras, to graded \textit{skew} Clifford algebras. In particular, in \cite{NVZ}, we show how GSCAs come into play with regard to most quadratic regular algebras of global dimension three and in \cite{VVe1, VVe2}, we review results on a notion of rank for a noncommutative quadratic form and discuss a generalization to Theorem \ref{vvw}. 

\subsection{Graded Skew Clifford Algebras}\label{}
\medskip

For $\{i,\ j\} \subset \{1, \ldots , n\}$, let $\mu_{ij} \in \kk^{\times}$ such that $\mu_{ij}\mu_{ji} = 1$ for all $i, j$ where $i \ne j$. We write $\mu = (\mu_{ij}) \in M(n,\ \kk)$ and $S$ as in \cite{CaV} for the quadratic $\kk$-algebra on degree-one generators $z_1, \ldots, z_n$ and defining relations $z_j z_i = \mu_{ij} z_i z_j$ for all $i$, $j = 1, \ldots, n$, where $\mu_{ii} = 1$ for all $i$. That is, \[S = \frac{T(V)}{\langle z_jz_i - \mu_{ij}z_iz_j : i, j = 1, \ldots, n\rangle},\] where $T(V)$ is the tensor algebra on $V$ = span$\{z_1, \ldots, z_n\}$. We set $U \subset S_1 \otimes_{\Bbbk} S_1$ to be the span of the defining relations of $S$ and write $z = (z_1, \ldots, z_n)^T$.

\begin{defn}\cite[\S1.2]{CaV}\label{MuSymm}
\begin{enumerate}
\item[(a)]
With $\mu$ and $S$ as above, a (noncommutative) quadratic form is any element of $S_2$.  Identifying $\mathbb{P}(S_1^{\ast})$ with $\mathbb{P}^{n-1}$, we note that the subscheme of the zero locus $Z \subset \mathbb{P}^{n-1} \times \mathbb{P}^{n-1}$ of $U$ consisting of all points in $Z$ on which a quadratic form $Q$ vanishes is called the quadric determined by $Q$.
\item[(b)] A matrix $M \in M(n,\ \kk)$ is called $\mu$-symmetric if $M_{ij} = 
\mu_{ij}M_{ji}$ for all $i$, $j = 1, \ldots , n$. 
\end{enumerate}
\end{defn}
We write $M^{\mu}(n, \kk)$ for the set of $\mu$-symmetric matrices in $M(n, \kk)$ and note that if $\mu_{ij} = 1$ 
for all $i, j$, then $M^{\mu}(n, \kk)$ is the set of all symmetric matrices.

\medskip


\begin{defn}\cite{CaV}\label{GSCA}
A {\em graded skew Clifford algebra} $A = A(\mu, M_1, \ldots , M_n)$ 
associated to $\mu$ and $M_1$, $\ldots ,$ $M_n \in M^{\mu}(n,\ \kk)$ is a $\mathbb{Z}$-graded $\kk$-algebra
on degree-one generators $x_1, \ldots , x_n$ and on degree-two generators
$y_1, \ldots , y_n$ with defining relations given by:
\begin{enumerate}
\item[(a)] $\displaystyle x_i x_j + \mu_{ij} x_j x_i = \sum_{k=1}^n (M_k)_{ij} y_k$
           for all $i, j = 1, \ldots , n$, and
\item[(b)] the existence of a normalizing sequence 
           $\{ y_1', \ldots , y_n'\}$ that spans 
	   $\kk y_1 + \cdots + \kk y_n$.
\end{enumerate}
\end{defn}

\begin{rem}\label{linIndepMuSymms}
If $A$ is a GSCA, then \cite[Lemma 1.13]{CaV} implies that $y_i \in (A_1)^2$ for all $i = 1, \ldots, n$ if and only if $M_1, \ldots, M_n$ are linearly independent. Henceforth, we assume that $M_1, \ldots, M_n$ are linearly independent.
\end{rem}

\begin{ex}
Consider the following four $\mu$-symmetric matrices: 
$M_1 = \left[\begin{smallmatrix}
0 & 1 & 0 & 0\\
\mu_{21} & 0 & 0 & 0\\
0 & 0 & 0 & 0\\
0 & 0 & 0 & 2
\end{smallmatrix} \right]$,
$M_2 = \left[\begin{smallmatrix}
0 & 0 & 1 & 0\\
0 & 2 & 0 & 0\\
\mu_{31} & 0 & 0 & 0\\
0 & 0 & 0 & 0
\end{smallmatrix} \right]$,
$M_3 = \left[\begin{smallmatrix}
0 & 0 & 0 & 1\\
0 & 0 & 0 & 0\\
0 & 0 & 2 & 0\\
\mu_{41} & 0 & 0 & 0
\end{smallmatrix} \right]$,
$M_4 = \left[\begin{smallmatrix}
2 & 0 & 0 & 0\\
0 & 0 & 1 & 0\\
0 & \mu_{32} & 0 & 0\\
0 & 0 & 0 & 0
\end{smallmatrix} \right]$.
\end{ex} \noindent Using Definition \ref{GSCA}(a), we obtain the following relations:
\[
\begin{array}{llr}
x_1 x_2 + \mu_{12} x_2 x_1 = x_4^2, \qquad & x_1 x_4 + \mu_{14} x_4 x_1 = x_2^2 , \qquad & 
x_2 x_4 = -\mu_{24}x_4 x_2,\\[3mm]
x_1x_3 + \mu_{13}x_3 x_1 = x_3^2, & \quad x_2 x_3 + \mu_{23}x_3 x_2= x_1^2, & 
x_3x_4 = -\mu_{34}x_4x_3.\\[3mm]
\end{array}\] 
We know that, $A'$, the $\Bbbk$-algebra on $x_1, x_2, x_3, x_4$ factored out by the ideal generated by the above relations, is mapped onto a GSCA. However, since $x_1^2$, $x_2^2$, $x_3^2$, and $x_4^2$ are normal in $A'$, then $A'$ is in fact a GSCA.  

\medskip

\begin{defn}\label{quadric}\hfill\\
\indent (a) \cite{CaV} \  
The span of quadratic forms $Q_1, \ldots , Q_m\in S_2$ will be called the 
{\em quadric system} associated to $Q_1, \ldots , Q_m$. If a quadric system is given by a normalizing sequence in $S$, then it is
called a {\em normalizing quadric system}.\\
\indent (b) \cite{CaV} \  
We define a {\em left base point} of a quadric system $\mathfrak{Q} \subset S_2$ to be 
any left base-point module over $S/ \langle \mathfrak{Q} \rangle$;
that is, to be any 1-critical graded left module over $S/\langle \mathfrak{Q} \rangle$
that is generated by its homogeneous degree-zero elements and which
has Hilbert series $H(t) = c/(1-t)$, for some $c \in \mathbb{N}$.
We say a quadric system is {\em left base-point free} if it has no left 
base points. Respectively, the same is true for right base point and right base-point free.\\
\indent (c) \cite{ATV2} In (b), if $c = 1$, then a left base-point module is cyclic and is a (left) point module.\\
\indent (d) \cite{ATV1} In (b), if $c \ge 2$, then a left base-point module is a (left) fat point module.\\
\end{defn}

\begin{rem}
\begin{enumerate}[(a)]
\item The quadric system, $\{q_1, \ldots, q_n\}$, in Theorem \ref{BPFToRegGSCAs} below can be obtained by associating the algebra $S$, defined at the beginning of this section, to the GSCA via $\mu$. The isomorphism $M^{\mu}(n, \Bbbk) \to S_2$ defined by $M \mapsto (z_1, \cdots, z_n)M(z_1, \cdots, z_n)^T$ associates the quadric system $\{q_1, \ldots, q_n\}$ to the $\mu$-symmetric matrices $M_1, \cdots, M_n$.
\item In \cite[Corollary 11]{CaV}, a normalizing quadric system $\mathfrak{Q}$ is right base-point free if and only if $\dim_{\Bbbk}(S/ \langle \mathfrak{Q} \rangle) < \infty$. Thus, a normalizing quadric system $\mathfrak{Q}$ is right base-point free if and only if it is left base-point free. 
\item For an understanding of 1-critical graded modules and other similar concepts, see Levasseur's paper \cite{Le}.  
\end{enumerate}
\end{rem}

\medskip

\begin{thm}\cite[Theorem 4.2]{CaV}\label{BPFToRegGSCAs}
A graded skew Clifford algebra $A = A(\mu, M_1, \ldots, M_n)$ is a quadratic Auslander-regular algebra of global dimension $n$ that satisfies the Cohen-Macaulay property with Hilbert series $1/(1-t)^n$ if and only if the quadric system associated to $\{q_1, \ldots, q_n\}$ is normalizing and base point free; in this case, $A$ is an Artin-Schelter regular $\mathbb{N}$-graded $\Bbbk$-algebra of global dimension $n$, a noetherian domain and unique up to isomorphism.
\end{thm}

\medskip

\begin{rem}\label{W}
\begin{enumerate}[(a)]
\item Hereafter, we assume that the quadric system associated to $q_1, \ldots, q_n$ is normalizing and base-point free. Theorem \ref{BPFToRegGSCAs} allows us to write $A = T(V)/\langle W \rangle$, where $V = S_1^{\ast}$ and $W \subset T(V)_2$. We write the Koszul dual \cite{ST} of $A$ as $T(S_1)/\langle W^{\perp} \rangle = S/\langle q_1, \ldots, q_n \rangle$. We note that $\{x_1, \ldots, x_n\}$ is the dual basis in $V$ to the basis $\{z_1, \ldots, z_n\}$ of $S_1$ and we write $\sum_{i, j} \alpha_{ijm} (x_ix_j + \mu_{ij}x_jx_i) = 0$ for the defining relations of $A$, where $\alpha_{ijm} \in \kk$ for all $i, j, m$, and $1 \le m \le n(n-1)/2$. 
\item By \cite[Lemma 5.1]{CaV}, there is a one-to-one correspondence between the set of pure tensors in $\mathbb{P}(W^{\perp})$, that is, $\{a \otimes b \in \mathbb{P}(W^{\perp}): a, b \in S_1\}$ and the zero locus $\Gamma$ of $W$ given by $\Gamma = \{(a,\, b) \in \mathbb{P}(S_1) \times \mathbb{P}(S_1) : w(a,\, b) = 0 \mbox{ for all } w \in W\}$. 
\item If a quadratic $\Bbbk$-algebra is Auslander-regular and has polynomial growth, then it is Artin-Schelter regular. (See \cite{Le, LeS} and appendix \ref{appendix} in this paper.)
\end{enumerate}
\end{rem}


\medskip

\subsection{Classification of quadratic regular algebras of global dimension three via GSCAs}


We now survey, Nafari, Vancliff, and Zhang's results \cite{NVZ} on the classification of most quadratic regular algebras of global dimension three as graded skew Clifford algebras. The authors consider all possible `types' of quadratic regular algebras of global dimension three and show they can be viewed as regular quadratic graded skew Clifford algebras. The analysis is split into two main cases: quadratic regular algebras of global dimension three with reducible or non-reduced point scheme (point scheme contains a line) and quadratic regular algebras of global dimension three with irreducible reduced point scheme. Recall that the point scheme of a quadratic regular algebra represents the functor of point modules. By \cite{ATV1}, the point scheme of a quadratic regular algebra of global dimension three is a cubic divisor in $\mathbb{P}^2$ or is all of $\mathbb{P}^2$. The next result relates algebras whose point scheme contains a line, that is whose point scheme is reducible or non-reduced, to GSCAs.

\begin{thm}\cite[Theorem 1.10]{NVZ}
Suppose $A$ is a quadratic AS-regular algebra of global dimension three, suppose $C \subset \mathbb{P}^2$ is its associated point scheme and suppose char$(\Bbbk) \ne 2$. 
If the point scheme of $A$ is reducible or non-reduced, then either
\begin{enumerate}
\item[(a)] $A$ is a twist, by an automorphism, of a GSCA, or
\item[(b)] $A$ is a twist, by a twisting system, of an Ore extension of a regular GSCA of global dimension two. 
\end{enumerate}
\end{thm}

This result provides us with several examples of quadratic AS-regular algebras of global dimension three via GSCAs. The following example \cite[Example 1.9]{NVZ} is that of a regular GSCA of global dimension three with a reducible or non-reduced point scheme, that is, a point scheme that contains a line.

\begin{ex}
Suppose that char$(\Bbbk) \ne 2$. Let 
$M_1 = \left[\begin{smallmatrix}
2 & 0 & 0\\
0 & 0 & 0\\
0 & 0 & 0
\end{smallmatrix} \right]$,
$M_2 = \left[\begin{smallmatrix}
0 & 0 & 0\\
0 & 2 & 0\\
0 & 0 & 0
\end{smallmatrix} \right]$, and
$M_3 = \left[\begin{smallmatrix}
0 & 1 & 0\\
\mu_{21} & 0 & 0\\
0 & 0 & 2
\end{smallmatrix} \right]$ and fix $\mu_{13}\mu_{23} = 1$ and $\mu_{13} + \mu_{12}\mu_{23} \ne 0$. Using Theorem~\ref{BPFToRegGSCAs}, observe that the $\mu$-symmetric matrices $M_1, M_2$ and $M_3$ yield a regular GSCA, $A$, of global dimension three on generators $x_1, x_2, x_3$ with defining relations 
\[x_1x_2 + \mu_{12}x_2x_1 = 0.\] 
\[x_1x_3 + \mu_{13}x_3x_1 = 0.\]
\[x_2x_3 + \mu_{23}x_3x_2 = x_3^2.\]  

\noindent The point scheme of $A$ is isomorphic to a subscheme $\mathcal{P}$ of $\mathbb{P}^2$. The subscheme $\mathcal{P}$ can computed by considering all the points $((a_1, a_2, a_3), (b_1, b_2, b_3))$ in $\mathbb{P}^2 \times \mathbb{P}^2$ such that \[(x_1x_2 + \mu_{12}x_2x_1)((a_1, a_2, a_3), (b_1, b_2, b_3)) = 0, \mbox{ and }\] \[(x_1x_3 + \mu_{13}x_3x_1)((a_1, a_2, a_3), (b_1, b_2, b_3)) = 0, \mbox{ and }\] \[(x_2x_3 + \mu_{23}x_3x_2 - x_3^2)((a_1, a_2, a_3), (b_1, b_2, b_3)) = 0.\] That is, solve
\[\left[\begin{smallmatrix}
\mu_{12}a_2 & a_1 & -a_3\\
\mu_{13}a_3 & 0 & a_1\\
0 & \mu_{23}a_3 & a_2
\end{smallmatrix} \right]
\left[\begin{smallmatrix}
b_1\\
b_2\\
b_3
\end{smallmatrix} \right] = 
\left[\begin{smallmatrix}
0\\
0\\
0
\end{smallmatrix} \right].\]

\medskip

To find solutions in $\mathbb{P}^2 \times \mathbb{P}^2$, we may assume that the determinant of the above $3 \times 3$ matrix is zero. This yields $((\mu_{13} + \mu_{12}\mu_{23})a_1a_2 + a_3^2)a_3 = 0$. That is, $\mathcal{P} = \mathcal{V}(((\mu_{13} + \mu_{12}\mu_{23})x_1x_2 + x_3^2)x_3)$, which is the union of a nondegenerate conic and a line.
\end{ex}


For the rest of this section, Nafari et al. consider regular algebras of global dimension three with a point scheme given by a cubic divisor $C$ that is reduced and irreducible. Such point schemes do not contain lines and thus, the cubic divisor has at most one singular point. Moreover, $x, y$ and $z$ are used for homogeneous degree-one linearly independent (commutative) coordinates in $\mathbb{P}^2$. It should be noted that algebras with nodal cubic and cuspidal cubic curves as point schemes are not discussed in \cite{ATV1} since they are not generic. 

In \cite[Lemma 2.1]{NVZ}, an irreducible cubic divisor, $C$ in $\mathbb{P}^2$, with a unique singular point is given by $\mathcal{V}(f)$ where either (a) $f = x^3 + y^3 + xyz$, or (b) $f = y^3 + x^2z$, or (c) $f = y^3 + x^2z + xy^2$ and if char$(\Bbbk) \ne 3$, then $f$ is given by (a) or (b). If $C$ is given by (a), then $C$ is called a nodal cubic curve and called a cuspidal cubic curve, otherwise. In both cases where the algebras' point scheme is given by either a nodal cubic curve or cuspidal cubic curve, the authors show that the algebras are in fact isomorphic to graded skew Clifford algebras or Ore extensions of regular GSCAs.

In \cite{AS} regular algebras of global dimension three with point scheme an elliptic curve are classified into four types: $A, B, E$ and $H$. These types along with their relations are summarized in the table below. The reader should note that some members of each type might not have an elliptic curve as their point scheme, but a generic member does. Nafari et al. \cite{NVZ} show that algebras of types $H$, $B$, and some of type $A$ can be related to graded \textit{skew} Clifford algebras and graded Clifford algebras. Algebras of type $E$ and an open subset of those of type $A$ were not tied to GSCAs.

\begin{center}
\begin{tabular}{ | m{1cm} | m{10em} | m{12em} | }
\hline Type & Defining Relations & Conditions \\
\hline \multirow{4}{7em}{H} & {$y^2=x^2$},& char$(\Bbbk) \ne 2$ \\
                                             & {$zy=-iyz$},&  \\
                                             & {$yx-xy=iz^2$},& \\
                                             & where $i$ is a primitive fourth root of unity.& \\
\hline \multirow{4}{7em}{B} & $xy+yx=z^2-y^2$,&$a(a - 1) \ne 0$, \\
                                            & $xy+yx=az^2-x^2$,& char$(\Bbbk) \notin \{2, 3\}$\\
                                            & $zx-xz=a(yz-zy)$,&\\
                                            & $a \in \Bbbk$, &\\
\hline \multirow{4}{7em}{A} & $axy + byx + cz^2 = 0$,& $abc \ne 0$, char$(\Bbbk) \ne 3$\\
					  & $ayz + bzy + cx^2 = 0$,& $(3abc)^3 \ne (a^3 + b^3 + c^3)^3$,\\
					  & $azx + bxz + cy^2 = 0$,& $a^3 = b^3 \ne c^3$, $b^3 = c^3 \ne a^3$, \\
					  & $a, b, c \in \Bbbk$.& $a^3 = c^3 \ne b^3, a^3 \ne b^3 \ne c^3$\\ 	      
\hline \multirow{4}{7em}{E} & $\gamma zx + xz = -\gamma^5 y^2$, &\\
					  & $yx + \gamma^4 xy = -\gamma^2z^2$,&\\
					  & $zy + \gamma^7 yz = -\gamma^8 x^2$,&\\
					  & where $\gamma$ is a primitive ninth root of unity.&\\
\hline	                                   
\end{tabular}
\end{center}

\bigskip

\noindent The algebras of type $H$, $B$, and $A$ will be denoted by $H$, $B$, and $A$, respectively. 

\begin{lem}\cite{NVZ}
\begin{enumerate}[(a)]
\item If char$(\Bbbk) \ne 2$, then the algebra $H$ is a regular graded skew Clifford algebra.
\item If char$(\Bbbk) \ne 2$, then the algebra $B$ is regular if and only if $a^2 - a + 1 \ne 0$; in this case, $B$ is a graded skew Clifford algebra and a twist of a Clifford algebra by an automorphism.
\item Suppose char$(\Bbbk) \ne 2$. If $a^3 = b^3 \ne c^3$, then $A'$ is a regular graded skew Clifford algebra and a twist of a graded Clifford algebra by an automorphism. If $b^3 = c^3 \ne a^3$ or if $a^3 = c^3 \ne b^3$, then $A'$ is a twist of a regular graded skew Clifford algebra by an automorphism.
\end{enumerate}
\begin{proof} (Sketch of proof) 
The main idea for the proof to parts (a), (b), and (c) is to use Remark \ref{KDGCA} parts (i) and (ii) to compute the Koszul dual of $H$, $B$, and $A'$ respectively, and apply Theorem \ref{BPFToRegGSCAs} to a certain normalizing sequence with empty zero locus in $\mathbb{P}^2 \times \mathbb{P}^2$. The proof to part (c) uses a similar argument for each of the three cases except that the three cases reduce to one case when $A'$ is twisted.
\end{proof}
\end{lem}

\medskip

\begin{rem}
We note that that a GCA (discussed in \S \ref{sec2}) is a finite module over some commutative subalgebra $C$, while a GSCA is a finite module over a (likely noncommutative) subalgebra $R$. Nafari et al. \cite{NV} show that if a regular GSCA is a twist by an automorphism of a GCA, then the subalgebra $R$ is a skew polynomial ring and is a twist of the commutative subalgebra $C$ by an automorphism.     
\end{rem}

\subsection{Point Modules over GSCAs}\label{PtModsGSCAs}


In this section, we discuss results \cite{VVe1, VVe2} that connect point modules to GSCAs. We see in Theorem \ref{bigOne} using the definition of $\mu$-rank given below that the quadratic forms of $\mu$-rank at most two control point modules over GSCAs. This result generalizes Theorem \ref{vvw} that related (commutative) quadratic forms of rank at most two to point modules over GCAs. 
\medskip

As in previous sections, for $\{i,\ j\} \subset \{1, \ldots , n\}$, let $\mu_{ij} \in \kk^{\times}$ such that $\mu_{ij}\mu_{ji} = 1$ for all $i \ne j$ and we write $\mu = (\mu_{ij}) \in M(n,\ \kk)$ and \[S = \frac{T(V)}{\langle z_jz_i - \mu_{ij}z_iz_j : i, j = 1, \cdots, n\rangle},\] where $T(V)$ is the tensor algebra on $V$ = span$\{z_1, \cdots, z_n\}$ and we write $z = (z_1, \dots, z_n)^T$.

\medskip

We begin with a generalization of the notion of rank, called $\mu$-rank, for noncommutative quadratic forms.

\medskip

\begin{defn}\cite{VVe1}\label{defnDi}
Let $M =
\left[\begin{smallmatrix}
    a & d & e \\[1mm]
    \mu_{21}d & b & f \\[1mm]
    \mu_{31}e & \mu_{32}f & c \end{smallmatrix}
\right] \in M^{\mu}(3, \kk)$ and, for $1 \le i \le 8$,
define the functions $D_i: M^{\mu}(3, \kk) \to \kk$ by\\[-3mm]
\begin{gather*}
\begin{array}{ll}
D_1(M) = 4d^2 - (1 + \mu_{12})^2 ab, \qquad &
D_4(M) = 2(1 + \mu_{23})de - (1 + \mu_{12})(1 + \mu_{13})af,\\[2mm]
D_2(M) = 4e^2 - (1 + \mu_{13})^2 ac, &
D_5(M) = 2(1 + \mu_{12})ef - (1 + \mu_{13})(1 + \mu_{23})cd,\\[2mm]
D_3(M) = 4f^2 - (1 + \mu_{23})^2 bc, &
D_6(M) = 2(1 + \mu_{13})df - (1 + \mu_{12})(1 + \mu_{23})be,
\end{array}\\[2mm]
\begin{array}{l}
D_7(M) = (\mu_{23}cd^2 - 2def + be^2)
            (\mu_{13}\mu_{21}cd^2 - 2def + \mu_{12}\mu_{23}\mu_{31}be^2),
	         \\[3mm]
D_8(M) = \mu_{21}(d + X)(e - Y) + \mu_{23}\mu_{31} (d-X)(e+Y)-2af,
\end{array}
\end{gather*}
where $X^2 = d^2 - \mu_{12} ab$ \ and \ $Y^2 = e^2 - \mu_{13} ac$.
We call $D_1, \ldots , D_6$ the $2\times 2$ $\mu$-minors of $M$. Since the functions $D_7$ and $D_8$ will play a role analogous to the determinant of $M$, they will be called $\mu$-determinants of $M$.
\end{defn}

\begin{thm}\cite{VVe1}\label{muRk3VarsThm}
Let $Q = az_1^2 + bz_2^2 + cz_3^2 + 2dz_1z_2 + 2ez_1z_3 + 2fz_2z_3 \in
S_2$, where $a, \ldots, f \in \kk$, and let $M \in M^{\mu}(3,\ \kk)$ be the 
$\mu$-symmetric matrix associated to $Q$. 
\begin{enumerate}[(a)]
\item[{\rm (a)}] 
There exists $L \in S_1$ such that $Q = L^2$
if and only if $D_i(M) = 0$ for all $i = 1, \ldots , 6$.
\item[{\rm (b)}] 
\begin{enumerate}
\item[{\rm (i)}] 
If $a = 0$, then there exists $L_1$, $L_2 \in S_1$ such that
$Q = L_1 L_2$ if and only if $D_7(M) = 0$;
\item[{\rm (ii)}] 
if $a \neq 0$, then there exists $L_1$, $L_2 \in S_1$ such that
$Q = L_1 L_2$ if and only if $D_8(M) = 0$ for some $X$ and $Y$ satisfying
$X^2 = d^2 - \mu_{12} ab$ \ and \ $Y^2 = e^2 - \mu_{13} ac$.
\end{enumerate}
\end{enumerate}
\end{thm}

\begin{defn}\cite{VVe1}\label{muRank3Vars}
Let $Q = az_1^2 + bz_2^2 + cz_3^2 + 2dz_1z_2 + 2ez_1z_3 + 2fz_2z_3 \in
S_2$, where $a, \ldots , f \in \kk$, with $a = 0$ or 1, 
let $M \in M^{\mu}(3,\ \kk)$ be the $\mu$-symmetric matrix associated to $Q$ 
and let $D_i : M^{\mu}(3,\ \kk) \to \kk$, for $i = 1, \ldots , 8$, be 
defined as in Definition~\ref{defnDi}.
If $n=3$, we define the function $\mu$-rank $: S_2 \to \mathbb{N}$ as follows:
\begin{enumerate}[(a)]
  \item if $Q = 0$, we define $\mu$-rank$(Q) = 0$;
  \item if $Q \ne 0$ and if $D_i(M) = 0$ for all $i = 1, \ldots , 6$, 
  we define $\mu$-rank$(Q) = 1$;
  \item if $D_i(M) \ne 0$ for some $i = 1, \ldots , 6$ and if 
          \[(1 - a)D_7(M) + aD_8(M) = 0,\] 
	  we define $\mu$-rank$(Q) = 2$;
  \item if $(1 - a)D_7(M) + aD_8(M) \neq 0$, we define $\mu$-rank$(Q) = 3$.
\end{enumerate}
\end{defn}

\medskip

In general, Definition \ref{muRank3Vars} suggests Definition \ref{muRankGen} for a notion of $\mu$-rank at most two for (noncommutative) quadratic forms on $n$ generators for any $n \in \mathbb{N}$.

\begin{defn}\cite{VVe2}\label{muRankGen}
Let $S$ be as defined at the beginning of \S\ref{PtModsGSCAs}, where $n$ is an arbitrary
positive integer, and let $Q \in S_2$.
\begin{enumerate}
\item[{\rm (a)}]
If $Q = 0$, we define $\mu$-rank$(Q) = 0$.
\item[{\rm (b)}]
If $Q = L^2$ for some $L \in S_1^\times$, we define $\mu$-rank$(Q) = 1$.
\item[{\rm (c)}]
If $Q \neq L^2$ for any $L \in S_1^\times$, but $Q = L_1 L_2$ where
$L_1$, $L_2 \in S_1^\times$, we define $\mu$-rank$(Q) = 2$.
\end{enumerate}
\end{defn}

\begin{rem}\label{tau}
Let $\tau: \mathbb{P}(M^{\mu}(n, \Bbbk)) \to \mathbb{P}(S_2)$ be defined by $\tau(M) = z^TMz$. Hereafter, we fix $M_1, \ldots, M_n \in M^{\mu}(n, \Bbbk)$. For each $k = 1, \ldots, n$, we fix representatives $q_k = \tau(M_k)$. Moreover, if $M$ is a $\mu$-symmetric matrix in $\mathbb{P}(M^{\mu}(n, \Bbbk))$ and if $\mu$-rank$(T(M)) \le 2$, then we define $\mu$-rank$(M)$ to be $\mu$-rank of $\tau(M)$.
\end{rem}

\begin{ex}
We consider the quadric system $\mathfrak{Q}$ given by $q_1 = z_1z_2, q_2 = z_3^2, q_3 = z_1^2 - z_2z_4$, and $q_4 = z_2^2 + z_4^2 - z_2z_3$ in \cite[\S 5.3]{CaV}. We observe that $\mathfrak{Q}$ is a normalizing and base-point free quadric system if and only if $\mu_{34}^2 = \mu_{23} = 1$ and $\mu_{34} = \mu_{24} = (\mu_{14})^2 = (\mu_{13})^2$. The quadric system is base-point free since the ideal $I$ it generates contains $z_1^3, z_2^5, z_3^2, z_4^5$. Thus, $\dim_{\Bbbk}(S/I) < \infty$ and \cite[Corollary 11]{CaV} implies that $\mathfrak{Q}$ is a base-point free quadric system.
\end{ex}

\medskip

In order to relate point modules to GSCAs, Vancliff and the author consider the quadratic forms of $\mu$-rank at most two that are in the Koszul dual (see Remark \ref{KDGCA}) of the GSCA. 

A point $(a, b)$ is associated to a $\mu$-symmetric matrix via the map $\Phi$ from $\mathbb{P}^{n-1} \times \mathbb{P}^{n-1}$ to $\mathbb{P}(M^{\mu}(n, \Bbbk))$ and the image of the map $\Phi$ restricted to $\Gamma$ (see Remark \ref{W}(b)) is shown to consist of $\mu$-symmetric matrices in $\mathbb{P}(\sum_{k=1}^n \Bbbk M_k)$ of $\mu$-rank at most two. Since a noncommutative quadratic form factors in at most two ways \cite[Theorem 2.8]{VVe2}, the following result follows. 

\begin{thm}\cite{VVe2}\label{bigOne}
If the quadric system $\{q_1, \ldots , q_n\}$ associated to the GSCA, $A$, is normalizing and base-point free, then the number of isomorphism classes of left (respectively, right) point modules over $A$ is equal to $2 f_2 + f_1 \in \mathbb{N} \cup \{0,\, \infty\}$, where $f_j$ denotes the number of matrices~$M$ in $\mathbb{P}\big(\sum_{k=1}^n \kk M_k\big)$ such that $\mu$-rank$(M) \le 2$ and such that $\tau(M)$ factors in $j$ distinct ways (up to a nonzero scalar multiple). 
\end{thm}

\begin{ex}
We consider a GSCA in \cite[\S5.3]{CaV} with 
$n = 4$, where\\[-3mm]
\begin{gather*} 
q_1 = z_1 z_2, \qquad q_2 = z_3^2, \qquad
q_3 = z_1^2 - z_2z_4, \qquad q_4 = z_2^2 + z_4^2 - z_2z_3,\\
\mu_{23} = 1 = -\mu_{34}, \qquad (\mu_{14})^2 = \mu_{24} = -1,\qquad
\mu_{13} = -\mu_{14},
\end{gather*} 
This quadric system is normalizing and base-point free. 
By Theorem \ref{BPFToRegGSCAs}, the corresponding GSCA, $A$, is quadratic and 
regular of global dimension four, and is the $\kk$-algebra on generators
$x_1, \ldots, x_4$ with defining relations:\\[2mm]
\[
\begin{array}{llr}
x_1 x_3 = \mu_{14} x_3 x_1, \qquad & x_3 x_4 = x_4 x_3, \qquad & 
x_2 x_3 + x_3 x_2= - x_4^2,\\[3mm]
x_1 x_4 = -\mu_{14} x_4 x_1, & \quad x_4^2 = x_2^2, & 
x_2 x_4 - x_4 x_2= - x_1^2,\\[3mm]
\end{array} 
\]
and has exactly five nonisomorphic point modules, two
of which correspond to $q_1 = z_1 z_2 = z_2 z_1$.
The other three point modules correspond to two quadratic forms in 
$\mathbb{P}\big(\sum_{k=1}^4 \kk q_k\big)$
that have $\mu$-rank one, namely
\[
q_2 = z_3^2 \qquad \text{and} \qquad 
q_2 + 4 q_4 = (z_2 -\frac{z_3}{2} + z_4)^2  = 
(-z_2 + \frac{z_3}{2} + z_4)^2,
\]
where the latter quadratic form clearly factors in two distinct ways.
Hence, $A$ has a finite number of point modules even though two distinct 
elements of $\mathbb{P}\big(\sum_{k=1}^4 \kk q_k\big)$ have $\mu$-rank one.
\end{ex}



\section{Current and Future Work} \label{sec6}

One of the goals for the work on GCAs and the introduction of GSCAs was to provide examples for quadratic regular algebras of global dimension three and four and this was the motivation for Section \ref{sec3}. Another goal is to provide candidates for generic quadratic regular algebras of global dimension four so as to contribute towards the classification of quadratic regular algebras of global dimension four. Indeed, two of the families of GSCAs discussed in \cite[\S 5]{CaV} are candidates for generic quadratic regular algebras of global dimension four. As discussed in \cite{V}, these algebras should have a point scheme with exactly twenty points and a one-dimensional line scheme. 

Using techniques that involve Pl\"{u}cker coordinates, \cite{ChV} analyzed the line scheme of this family of algebras. They found that the line scheme consists of one spatial elliptic curve, four planar elliptic curves, and a subscheme in a $\mathbb{P}^3$ consisting of the union of two nonsingular conics. In \cite{TV}, the authors analyze yet another family of algebras with a point scheme consisting of exactly twenty points and a one-dimensional line scheme and their analysis leads to a line scheme consisting of one spatial elliptic curve, one spatial rational curve, two planar elliptic curves, and two subschemes where each is the union of a nonsingular conic and a line. This leads the authors to conjecture that the line scheme of the most generic quadratic regular algebra of global dimension four should be isomorphic to the union of two spatial elliptic curves and four planar elliptic curves. 

\appendix
\section{Standard Definitions}\label{appendix}
\subsection{\textit{Global Dimension of a Graded Connected Algebra \cite{ATV1}}}
A graded connected (meaning $A_0 = \Bbbk$) algebra $A$ has global dimension $d$ if every left $A$-module and every right $A$-module has projective dimension at most $d$ and at least one left module and at least one right module has projective dimension equal to $d$.

\subsection{{\textit{Polynomial Growth of a Graded Algebra $A$ \cite{ATV1}}}}
Suppose $A$ is an $\mathbb{Z}$-graded algebra such that $A = \displaystyle\oplus_{i \ge 0} A_i$. The algebra $A$ has polynomial growth if $\dim_{\kk} A_n \le cn^{\delta}$ for some positive real numbers $c$ and $\delta$ for all $n$.

\subsection{\textit{Gorenstein} \cite{AS}}
An $\mathbb{Z}$-graded algebra $A$, of finite global dimension, is Gorenstein if:
\begin{enumerate}[(i)]
\item the projective modules appearing in a minimal resolution of the left trivial module $_{A}\Bbbk$ are finitely generated, and
\item the transposed complex (or the ``dual sequence") obtained by applying the functor $M \rightsquigarrow M^{\ast}$ = Hom$_A(M, A)$ to a minimal resolution of $_{A}\Bbbk$ is a resolution of a graded right module isomorphic to the right trivial module $\Bbbk_A$.  
\end{enumerate}

\subsection{{\textit{Hilbert series} of a module \cite{ATV2}}}
The Hilbert series of a graded $\mathbb{Z}$-module or a $\mathbb{Z}$-graded $\Bbbk$-vector space $M = \oplus M_n$ is the formal series
\[h_M(t) = \sum_n (\dim_{\Bbbk} M_n)t^n.\]

\subsection{{\textit{Normal Element }\cite{MR}}}
An element $a$ of a ring $R$ is a normal element if $aR = Ra$.

\subsection{{\textit{Normalizing Sequence} \cite{MR}} }
A sequence $a_1, ..., a_n$ of elements of a ring $R$ is called a normalizing sequence if:
\begin{enumerate}[(i)]
\item $a_{1}$ is a normal element of $R$,
\item for each $j \in \{1, ..., n-1\}$ the image of $a_{j+1}$ in $\frac{R}{\sum_{i=1}^j a_iR}$ is a normal element, and
\item $\sum_{i=1}^n a_iR \ne R$.
\end{enumerate}

\subsection{\textit{Auslander-Regular} Algebra \cite{Le}}
Let $A$ be a noetherian $\Bbbk$-algebra.
\begin{enumerate}[(i)]
\item An $A$-module $M$ satisfies the Auslander-condition if, for all $q \ge 0$, \[q \le \inf\{i: \mbox{Ext}_A^i(N, A) \ne 0\},\] for all $A$ submodules $N$ of Ext$^q(M, A)$.
\item The algebra $A$ is Auslander-regular of global dimension $d$ if gldim$(A) = d < \infty$ and every finitely generated $A$-module satisfies the Auslander-condition.
\end{enumerate}



\begin{thebibliography}{1}
%
%
\bibitem[AS]{AS}
{\sc M.~Artin and W.~.F.~Schelter}, Graded Algebras of Global Dimension 3, {\it Advances in Math} {\bf 66} (1987), 171-216.

\bibitem[ATV1]{ATV1} 
{\sc M.~Artin,  J.~Tate and M.~Van den Bergh},  Some Algebras
Associated to Automorphisms of Elliptic Curves, {\it The Grothendieck
Festschrift} {\bf 1}, 33-85,
Eds.\ P.\ Cartier et al., Birkh\"auser (Boston, 1990).

\bibitem[ATV2]{ATV2} 
{\sc M.~Artin, J.~Tate and M.~Van den Bergh}, Modules over Regular
Algebras of Dimension~3, {\it Invent.\ Math.} {\bf 106} (1991), 335-388.

\bibitem[AL]{AL}
{\sc M.~Aubry and J.-M. Lemaire}, Zero Divisors in Enveloping Algebras of Graded Lie Algebras,
{\it J. Pure \& Appl. Alg.} {\bf 38} (1985), 159-166.

\bibitem[CK]{CK}
{\sc Z.~Chen and Y.~Kang}, Generalized Clifford Theory for Graded Spaces, {\it J.\ Pure\ App.\ Algebra} {\bf 220} (2016), 647-665.

\bibitem[CaV]{CaV} 
{\sc T.\ Cassidy and M.\ Vancliff}, Generalizations of Graded Clifford Algebras and of Complete
Intersections, {\it J. Lond. Math. Soc.} {\bf 81} (2010), 91-112. (Corrigendum: {\bf 90} No. 2 (2014), 631-636.)

\bibitem[ChV]{ChV}
{\sc R.~Chandler and M.\ Vancliff}, The One-Dimensional Line Scheme of a Certain Family of Quantum $\mathbb{P}^3$s, 
{\it J. Algebra} {\bf 439} (2015), 316-333.


\bibitem[KL]{KL}
{\sc G.~Krause and T.~Lenagan}, Growth of algebras and Gelfand-Kirillov dimension, Revised edition, {\it Graduate Studies in Mathematics}, {\bf 22}, American Mathematical Society, Providence, RI, 2000.

\bibitem[Lb]{Lb}
{\sc L.~Le~Bruyn}, Central Singularities of Quantum Spaces, {\it J.~Algebra}, {\bf 177} (1995), 142-153. 

\bibitem[Le]{Le}
{\sc T.~Levasseur}, Some Properties of Non-Commutative Regular Graded Rings, {\it Glasgow Math. J.}, {\bf 34} (1992), 277-300.

\bibitem[LeS]{LeS}
{\sc T.~Levasseur and S.~P.~Smith}, Modules over the 4-Dimensional Sklyanin Algebra, {\it Bulletin de la S.M.F.}, {\bf 121} No.~1 (1993), 35-90.

\bibitem[L]{L}
{\sc P.~Lounesto}, Clifford Algebras and Spinors, {\it Lond.\ Math.\ Soc. Lect. Series}, {\bf 286} (2001).

\bibitem[MR]{MR}
{\sc J.~C.~McConnell and J.~C.~Robson}, Noncommutative Noetherian Rings, {\it Graduate Studies in Mathematics, American Mathematical Society}, {\bf 30} (2001).

\bibitem[NV]{NV}
{\sc M.~Nafari and M.~Vancliff}, Graded Skew Clifford Algebras that are Twists of Graded Clifford Algebras, {\em Comm. Alg.}
{\bf 43} No. 2 (2015), 719-725.  

\bibitem[NVZ]{NVZ}
{\sc M.~Nafari, M.~Vancliff and Jun Zhang}, Classifying Quadratic Quantum
$\mathbb{P}^2$s by using Graded Skew Clifford Algebras, {\it J.~Algebra}, {\bf 346}
No.~1 (2011), 152-164.

\bibitem[ST]{ST}
{\sc B.~Shelton and C.~Tingey}, On Koszul Algebras and a New Construction of Artin-Schelter Regular Algebras, {\em J.~Alg.}, 
{\bf 241} (2001), 789-798.

\bibitem[SV]{SV}
{\sc D. R.~Stephenson and M.~Vancliff}, Constructing Clifford Quantum $\mathbb{P}^3$s with Finitely Many Points, {\it J. Alg.},
{\bf 312} No. 1 (2007), 86-110.

\bibitem[TV]{TV}
{\sc D.~Tomlin and M.~Vancliff}, The One-Dimensional Line Scheme of a Family of Quadratic Quantum $\mathbb{P}^3$'s, {\it Work in progress}.

\bibitem[V]{V}
{\sc M.~Vancliff}, The Interplay of Algebra and Geometry in the Setting of Regular Algebras, in ``Commutative Algebra and Noncommutative Algebraic Geometry,''
{\it MSRI Publications}, {\bf 67} (2015), 371-390.

\bibitem[VVr]{VVr}
{\sc M.~Vancliff and K.~Van Rompay}, Four-dimensional Regular Algebras with Point Scheme a Nonsingular Quadric in $\mathbb{P}^3$, {\it Comm. Alg.}, 
{\bf 28} No. 5 (2000), 2211-2242.

\bibitem[VVe1]{VVe1}{\sc M.~Vancliff and P.~P.~Veerapen}, Generalizing the 
Notion of Rank to Noncommutative Quadratic Forms, {\em in}
``Noncommutative Birational Geometry, Representations and
Combinatorics,'' Eds.\ A.~Berenstein and V.~Retakh, {\em Contemporary 
Math.}\ {\bf 592} (2013), 241-250.

\bibitem[VVe2]{VVe2}
{\sc M.\ Vancliff and P.\ P.\ Veerapen}, Point Modules over Regular Graded Clifford Algebras, {\em J.~Algebra}, 
{\bf 420} (2014), 54-64. 

\bibitem[VVW]{VVW}
{\sc M.~Vancliff, K.~Van Rompay and L.~Willaert}, Some Quantum $\mathbb{P}^3$s
with Finitely Many Points, {\it Comm.\ Alg.} {\bf 26} No.~4 (1998), 1193-1208.

\end{thebibliography}
\end{document}